\documentclass[aap,MSNbibl,dvips]{arximspdf}

\doi{10.1214/13-AAP989} 
\volume{25}
\issue{1}
\pubyear{2015}
\firstpage{104}
\lastpage{115}

\makeatletter
\newcommand{\rrvert}{\vert}
\newcommand{\llvert}{\vert}
\newtheorem{theorem}{Theorem}
\newtheorem{lemma}{Lemma}
\makeatother

\begin{document}
\begin{frontmatter}

\title{Spatial Moran models\\
I. Stochastic tunneling in the neutral case}
\runtitle{Tunneling in spatial Moran models}

\begin{aug}
\author[A]{\fnms{Richard} \snm{Durrett}\corref{}\thanksref{T1}\ead[label=e1]{rtd@math.duke.edu}}
\and
\author[A]{\fnms{Stephen} \snm{Moseley}\thanksref{T2}}
\runauthor{R. Durrett and S. Moseley}
\affiliation{Duke University}
\address[A]{Department of Mathematics\\
Duke University\\
Box 90320\\
Durham, North Carolina 27708-0320\\
USA\\
\printead{e1}} 
\end{aug}
\thankstext{T1}{Supported in part by NIH Grant 5R01GM096190.}
\thankstext{T2}{This work was part of his 2011 Ph.D. thesis in Applied
Mathematics at Cornell University.}

\received{\smonth{12} \syear{2012}}
\revised{\smonth{11} \syear{2013}}

%
\begin{abstract}
We consider a multistage cancer model in which cells are arranged in a
$d$-dimensional integer lattice.
Starting with all wild-type cells, we prove results about the
distribution of the first time when two
neutral mutations have accumulated in some cell in dimensions $d \ge
2$, extending work done by Komarova [\textit{Genetics} \textbf{166} (2004) 1571--1579] for $d=1$.
\end{abstract}

%

\end{frontmatter}

\section{Introduction}\label{sec1}

The accumulation of mutations is important not only for cancer
initiation, progression, and metastasis, but also for the emergence of
acquired resistance against chemotherapeutics, radiation therapy, or
targeted drugs. For this reason there is a large and growing literature
on the waiting time $\tau_k$ until some cell has acquired $k$
prespecified mutations. In all the models we consider, type $i$
individuals mutate to type $(i+1)$ at rate $u_{i+1}$. The dynamics
considered have most often been studied in multi-type Moran models with
a homogeneously mixing population of constant size. Here we will
concentrate on how results change when one considers a spatial Moran
models, and as is the case for much earlier work we will concentrate on
the behavior of $\tau_2$.

We suppose that cells of type 0 and type 1 have relative fitness 1 and
$\lambda$. Since we will only consider the waiting time for the first
type 2, the relative fitness of type 2's is not important. In this work
we will consider situation in which $\lambda$ is so close to 1 that
the mutations are essentially neutral. For cancer applications, this is
a restrictive assumption, and it will be removed in the companion paper
(part II) by Durrett, Foo and Leder \cite{DFL}. However, the current
result applies to the important case of tumor suppressor genes. In that
case, when both copies of the gene are inactivated trouble develops,
but while there is one working copy the cell can function normally.

We begin by recalling results for the Moran model in a homogeneously
mixing population of size $N$.
Here and in what follows the mutation rates $u_i$ and selection
coefficient $\lambda$ depend on $N$,
even though this is not indicated in the notation, and
we write $a_N \ll b_N$ if $a_N/b_N \to0$ as $N\to\infty$. The next
result made its first appearance
on page~16,230 of Nowak et al. \cite{N02}. Since then it has appeared in
print a number of times:
\cite{KSN,N04,IMN,IMKN}, and in Nowak's
excellent book \cite{N06} on Evolutionary Dynamics.

\begin{theorem} \label{wait2}
In the neutral case of the Moran model, $\lambda=1$, if we assume that
%
\begin{equation}
\frac{1}{\sqrt{u_2}} \ll N \ll\frac{1}{u_1} \label{th0cond}
\end{equation}
and let $u_1, u_2 \to0$ then we have
\[
P\bigl( \tau_2 > t/Nu_1 u_2^{1/2}
\bigr) \to\exp( - t ). %
\]
The same conclusion holds if $|\lambda-1| \ll u_2^{1/2}$.
\end{theorem}

Durrett and Schmidt \cite{DS} applied these ideas to study regulatory
sequence evolution and to expose flaws in Michael Behe's arguments for
intelligent design. Durrett, Schmidt and Schweinsberg \cite{DSS}, see
also Schweinsberg \cite{Sch}, generalized this result to cover $\tau_k$.

The conditions in the result may look mysterious but they can be
derived by simple reasoning.
Here and throughout the paper and $f(u) \sim g(u)$ means $f(u)/g(u) \to
1$ as $u \to0$. Suppose first that $\lambda=1$.
\begin{longlist}[(A3)]
\item[(A1)] If we start the Moran model with $k \ll N$ type 1's and the rest
type 0, then the 1's behave like a critical branching process. The time
needed for the 1's to die out is $O(k)$ and the number of type-1 births
before they die out is $O(k^2)$. Thus we expect the first type 2 to
occur in a type-1 family that reaches size $k_1=O(1/\sqrt{u_2})$, and
hence has $O(k_1^2) = O(1/u_2)$ births. The condition $1/\sqrt{u_2}
\ll N$ in (1) guarantees $k_1 \ll N$.

\item[(A2)] Since the voter model is a martingale, the probability a type-1
mutation creates a family that reaches size $1/\sqrt{u_2}$ is $\sqrt
{u_2}$. More to the point a simple computation (consider what happens
at the first jump) shows that the probability a type-1 family gives
rise to a type 2 before it dies out is $\sim\sqrt{u_2}$. Since
mutations to type 1 occur at times of a rate $Nu_1$ Poisson process and
with probability $\sim\sqrt{u_2}$ give rise to a type 2,
it follows that if $\rho_2$ is the birth time of the type-1 mutant
that first gives rise to a type 2 then
\[
P\bigl( \rho_2 > t/Nu_1 u_2^{1/2}
\bigr) \to\exp( - t ). %
\]
To complete the proof we need to show that $\tau_2 - \rho_2 \ll\rho
_2$, and for this we need the condition $N \ll1/u_1$ in (1).

\item[(A3)] By the discussion of (A1), the first mutation will occur in a
family that reaches a size $O(1/\sqrt{u_2})$. If $|\lambda-1| \ll
u_2^{1/2}$, then computations with Girsanov's formula show that (in the
limit $u_2\to0$) the behavior of the Moran model, while it is
$O(1/\sqrt{u_2})$, is indistinguishable from the case with no drift.
\end{longlist}

The assumption of a homogeneously mixing cell populations simplifies
calculations considerably, but is not realistic for solid tumors. For
this reason, Komarova~\cite{K06} considered a spatial model, which is
very similar to one introduced much earlier by Williams and Bjerknes
\cite{WB}. Due to work of Bramson and Griffeath \cite{BG1,BG2}, the
second model is known to probabilists as the biased voter model.

In the usual formulation of the biased voter model, each site on the
\mbox{$d$-}dimen\-sional integer lattice $\mathbb{Z}^d$ can be in state 0 or 1
indicating the presence of a cell with relative fitness 1 or $\lambda>
1$. Cells give birth at a rate equal to their fitness, and the
offspring replaces a nearest neighbor chosen at random. When $\lambda
=1$ this is the voter model which was introduced independently by
Clifford and Sudbury \cite{CS} and Holley and Liggett \cite{HL}. For
a summary of what is known see Liggett \cite{L99}.

In the biased voter model births drive the process. In Komarova's
version cells die at rate 1 and are replaced by a copy of a nearest
neighbor chosen with probability proportional to its fitness. A site
with $n_i$ neighbors in state $i$ makes
\[
\begin{tabular}{cc}
transitions & at rate \\
$0 \to1$ & $\lambda n_1/(\lambda n_1 + n_0)$ \\
$1 \to0$ & $n_0/(\lambda n_1 + n_0)$
\end{tabular}
\]

In $d=1$ if the set of sites in state 1 is an interval $[\ell,r]$ with
$\ell< r$ then any site that can change has
$n_1=n_0=1$ so Komarova's model is a time change of the biased voter
model. In $d\ge2$ this is not exactly true.
However, we are interested in values of $\lambda=1+s$ where $s=0.02$
or even less, so we expect the two models
to have very similar behavior. In any case, the difference between the
two models is much less than their difference from reality,
so we will choose to study the biased voter, whose duality with
branching coalescing random walk (to be described below)
gives us a powerful tool for doing computations.

Since we want a finite cell population we will restrict our process to
be a subset of $(-L/2,L/2]^d$. Komarova \cite{K06} uses ``Dirichlet
boundary conditions'', that is, she assumes her space is an interval with
no cells outside, but this is awkward because the set of type-1 cells
may reach one end of the interval and then no further changes happen at
that end. To avoid this, we will use periodic boundary conditions,
that is, we consider $( \mathbb{Z}\bmod L)^d$. The resulting toroidal
geometry is a little strange for studying cancer. However, using $(
\mathbb{Z}\bmod L)^d$ has the advantage that the space looks the same
seen from any point. Our results will show that for the parameter
values the first type 2 will arise when the radius of the set of sites
occupied by 1's is $\ll L$ so the boundary conditions do not matter.

Let\vspace*{1pt} $\xi_s^0$ be the set of cells equal to 1 in the voter model with
no mutations from 0 to 1 on $\mathbb{Z}^d$ starting from a single type 1 at 0. Let $|\xi^0_s|$ be the number of cells in $\xi_s^0$, and let
%
\begin{equation}
\nu_d = 1 - E \exp\biggl(-u_2\int_0^{T_0}\bigl|
\xi_s^0\bigr| \,ds \biggr). \label{nuddef}
\end{equation}
This quantity, which is defined for the voter model without mutation, calculates
the probability, which depends on the dimension $d$, that a mutation to
type 1 gives rise to a type 2 before\vspace*{1pt} its family dies out.
To see why this is true note that the integral $\int_0^{T_0}|\xi
_s^0| \,ds$ gives the total number of man-hours in the type-1 family,
and conditional on this\vspace*{1pt} the number of mutations that will occur is
Poisson with mean $u_2\int_0^{T_0}|\xi_s^0| \,ds$.\vspace*{1pt}

Since mutations to type 1 in a population of $N$ cells occur at rate
$Nu_1$ this suggests that
%
\begin{equation}
P(\tau_2 > t)\to\exp(- Nu_1\nu_d t ).
\label{t2Neut}
\end{equation}
As we will explain in a moment, there is a constant $\gamma_d$ so that
$\nu_d \sim\gamma_d h_d(u_2)$ as $u_2\to0$ where
%
\begin{equation}
h_d(u) = \cases{u^{1/3}, &\quad$d=1$,
\vspace*{2pt}\cr
u^{1/2}
\log^{1/2}(1/u), &\quad$d=2$,
\vspace*{2pt}\cr
u^{1/2}, &\quad$d
\ge3$.}\label{hd}
\end{equation}
To state the result we need one more definition:
%
\begin{equation}
g_d(u) = \cases{u^{1/3}, &\quad$d=1$,
\vspace*{2pt}\cr
\log^{-1/2}(1/u), &\quad$d=2$,
\vspace*{2pt}\cr
1, &\quad$d\ge3$.}\label{gd}
\end{equation}

\begin{theorem} \label{th1}
In the neutral case of the biased voter model, $\lambda= 1$, if we assume
%
\begin{equation}
\frac{1}{h_d(u_2)} \ll N \ll\frac{g_d(u_2)}{u_1}, \label{th1cond}
\end{equation}
then there are constants $\gamma_d$ given in (\ref{gam1}) and (\ref
{gamd}) so that as $u_1, u_2 \to0$
\[
P\bigl(\tau_2 > t/ Nu_1\gamma_d
h_d(u_2)\bigr) \to\exp(- t). %
\]
The same conclusion holds if $|\lambda-1| \ll h_d(u_2)$.
\end{theorem}

In $d=1$ this result was proved by Komarova \cite{K06}, see her
equation (62) and assumption (60), then change notation
$u_1 \to u$, $u_2 \to u_1$. See also her survey paper~\cite{K07}.
Note that when $d\ge3$ the order of magnitude of the waiting time and
the assumptions are the same as in Theorem \ref{wait2}. In $d=2$ there
are logarithmic corrections to the behavior in Theorem \ref{wait2}, so
only in the case of $d=1$ (which is relevant to cancer in the mammary
ducts) does space make a substantial change in the waiting time.

The reasons for the conditions in Theorem \ref{th1} are the same as in
Theorem~\ref{wait2}.
\begin{longlist}[(B3)]
\item[(B1)] We will see that the mutation to type 2 will occur in a type-1
family that reaches size $k=O(1/h_d(u_2))$. The left-hand side assumption in
(\ref{th1cond}) implies that $k \ll N$, so the type-2 mutant arises
before the 1's reach fixation.

\item[(B2)] Let $\rho_2$ be the time of the first type-1 mutation that begins
the family that eventually leads to a type 2. Since mutations to type 1
occur at rate $Nu_1$ and lead to a type 2 with probability $\nu_d$, it
is easy to see that
\[
P(\rho_2 > t) \to\exp(-Nu_1 \nu_d t)
\]
so to prove the result we need to show that with high probability $\tau
_2 - \rho_2 \ll\rho_2$. As the reader will see, this is guaranteed
by the right-hand side assumption in (\ref{th1cond}).

\item[(B3)] As in the discussion of Theorem \ref{wait2}, once we know that
the mutation to type 2 will occur in a type-1 family that reaches size
$k=O(1/h_d(u_2))$, it follows that if $|\lambda-1| \ll h_d(u_2)$ then
(in the limit $u_2\to0$) the behavior of the size of the biased voter
$|\xi^0_t|$ is, while it is $O(1/h_d(u_2))$, indistinguishable from
the case with no drift.
\end{longlist}

\section{The key to the proof}\label{sec2}

The size of the voter model, when $|\xi^0_t|>0$, is a time change of
symmetric simple random walk,
with jumps happening at two times the size of the boundary $|\partial
\xi^0_t|$, which is the\vspace*{1pt} number of
nearest neighbor pairs with $x \in\xi^0_t$ and $y\notin\xi^0_t$.
The one-dimensional case is easy because when $\xi^0_t \neq\varnothing$
the boundary $|\partial\xi^0_t|=2$.
The key to the study of the process in $d \ge2$ is the observation
that there are constants $\beta_d$ so that
%
\begin{equation}
\bigl|\partial\xi^0_t\bigr| \sim_p \cases{2d
\beta_d \bigl|\xi^0_t\bigr|, &\quad$d\ge3$,
\vspace*{3pt}\cr
4
\beta_2 \bigl|\xi^0_t\bigr|/\log\bigl(|
\xi^0_t|\bigr), &\quad$d=2$,}\label{bdry}
\end{equation}
where $|\partial\xi^0_t| \sim_p f(|\xi^0_t|)$ means that when $|\xi
^0_t|$ is large, $|\partial\xi^0_t|/f(|\xi^0_t|)$
is close to 1 with high probability.

The intuition behind this result is that the voter model is dual to a
collection of coalescing random walks, so in $d \ge3$ neighbors of
points in $\xi^0_t$ will be unoccupied with probability $\approx\beta
_d$, the probability two simple random walks started at 0 and
$e_1=(1,0,\ldots, 0)$ never hit. In dimension $d=2$, the recurrence of
random walks implies that when $|\xi^0_t|=k$ is large, most neighbors
of points in $\xi^0_t$ will be occupied, but due to the fat tail of
the recurrence time sites will be vacant with probability $\sim\beta
_2/\log k$, where $\beta_2=\pi$.

Before we try to explain why (\ref{bdry}) is true, we will list an
important consequence. Let $T_k$ be the first time $|\xi^0_t|=k$. Let
\[
a_n = \cases{n^2, &\quad$d=1$,
\vspace*{2pt}\cr
2n\log n, &
\quad$d=2$,
\vspace*{2pt}\cr
n, &\quad$d\ge3$.} %
\]

\begin{lemma}
Let $\xi^0_t$ be the unbiased voter model (i.e., $\lambda= 1$)
starting from a single occupied site.
%
\begin{equation}
\biggl( \frac{|\xi^0_{T_{n\varepsilon} + a_n t}|}{n} \bigg| T_{n\varepsilon} <
\infty\biggr)
\Rightarrow(Y_t | Y_0 = \varepsilon), \label{procconv}
\end{equation}
where $\Rightarrow$ indicates convergence in distribution of the
stochastic processes and the limit has
\[
dY_t = \cases{ \sqrt{2} \,dB_t, &\quad$d=1$,
\vspace*{2pt}\cr
\sqrt{2
\beta_d Y_t} \,dB_t, &\quad$d\ge2$,} %
\]
where $B_t$ is a one-dimensional Brownian motion.
In $d=1$ the process is stopped when it hits 0. In $d\ge2$, 0 is an
absorbing boundary so we do not need to stop the process.
\end{lemma}

In $d=1$ this result is trivial. If one accepts (\ref{bdry}) then
(\ref{procconv}) can be proved easily by computing infinitesimal
means and variances and using standard weak convergence results.
In $d\ge2$, (\ref{bdry}) and (\ref{procconv}) are almost
consequences of work of Cox, Durrett and Perkins \cite{CDP}.
They speed up time at rate $a_n$, scale space by $1/\sqrt{a_n}$, and
assign each point occupied in the voter model mass $1/n$
to define a measure-valued diffusion $X^n$ which they prove converges
to super-Brownian motion. See their Theorem~1.2. (Their scaling is a little
different in $d=2$ but this makes no difference to the limit.)

Let $V'_{n,s}(x)$ be the fraction of sites adjacent to $x$ in state 0
at time $s$
(with the prime indicating that we multiply this by $\log n$ in $d=2$,
see page 196).
A key step in the proof in \cite{CDP} is to show, see (I1) on page
202, that
for nice test functions $\phi$
%
\begin{equation}
E \biggl[ \biggl( \int_0^T
X^n_s\bigl( \bigl\{ V'_{n,s} -
\beta_d \bigr\} \phi^2\bigr) \,ds \biggr)^{2}
\biggr] \to0, \label{bdryL2}
\end{equation}
where $X^n_s(f)$ denote the integral of the function $f$ against the
measure $X^n_s$.
The result in (\ref{bdryL2}) shows that when we integrate in time and
average in space
(multiplying by a test function to localize the average) then (\ref
{bdry}) is true.

From the convergence of the measure valued diffusion $X^n$ to
super-\break Brownian motion, (\ref{procconv}) follows by considering the
total mass.
Earlier we said (\ref{procconv}) is almost a consequence of \cite
{CDP}, since they start their process
from an initial measure [i.e., $O(n)$ initial 1's] while consider a
single occupied site and condition on reaching~$n\varepsilon$.
However, this defect can be remedied by citing the work of Bramson,
Cox and LeGall \cite{BCL}, who have a result, Theorem~4 on page 1012
that implies (\ref{procconv}) in $d \ge2$.

The result in (\ref{procconv}) is enough for Section~\ref{sec3}, but for the
calculations in Section~\ref{sec4} we will need a version of (\ref{bdry}).
In that section we will compute under the assumption that if $|\xi^0_t|=k$
%
\begin{equation}
\bigl|\partial\xi^0_t\bigr| = \cases{ 2d\beta_d k, &
\quad$d\ge3$,
\vspace*{2pt}\cr
4\beta_2 k/\log k, &\quad $d=2$.} \label{bdryA}
\end{equation}
If one wants to give a rigorous proof of the estimates there, then
small values of $k$, can be treated with the inequalities
\[
C k^{1/d} \le\bigl|\partial\xi^0_t\bigr| \le2 dk,
\]
and one can control large values of $k$ using (\ref{bdryL2}) and
estimates such as (J1) and (J2) on page 208 of \cite{CDP}.
We will assume (\ref{bdryA}) in order to avoid getting bogged down in
technicalities.

\section{Proof, part I}\label{sec3}

Let $\nu_d^\varepsilon$ be the probability defined in (\ref{nuddef})
ignoring mutations to type 2 that occur before $T_{n\varepsilon}$.
The size of the voter model, $|\xi^0_t|$, is a martingale, so if we
let $P_1$ to denote the law of the voter model starting
from one occupied site $P_1( T_{n\varepsilon} < \infty) = 1/n\varepsilon
$. Applying (\ref{procconv}) now,
%
\begin{equation}
\nu_d^\varepsilon\sim\frac{1}{n\varepsilon} \cdot\biggl[ 1 -
E_\varepsilon\exp\biggl( - n a_n u_2 \int
_0^{T_0} Y_s \,ds \biggr) \biggr],\label{nuasym}
\end{equation}
where $T_0 = \min\{ t\dvtx  Y_t=0\}$, $E_\varepsilon$ is the expected value
for $(Y_t | Y_0 = \varepsilon)$. We have
\[
n a_n = \cases{ n^3, &\quad$d=1$,
\vspace*{2pt}\cr
2n^2
\log n, &\quad$d=2$,
\vspace*{2pt}\cr
n^2, &\quad$d \ge3$.} %
\]
So if we set $n = 1/h_d(u_2)$ then (\ref{hd}) implies $n a_n u_2 \to
1$ and using (\ref{nuasym}) gives
\[
\nu_d^\varepsilon\sim h_d(u_2) \cdot
\biggl[ \frac{ 1 - E_\varepsilon
\exp( - \int_0^{T_0} Y_s \,ds ) }{\varepsilon} \biggr]. %
\]
Thus the type-2 mutation will occur in a family that reaches sizes
$O(1/h_d(u_2))$, and we must assume $1/h_d(u_2) \ll N$.

If we ignore the time to reach size $1/h_d(u_2)$, the time needed to
generate the type-2 mutation is, by (\ref{procconv}), of order
\[
a\bigl(1/h_d(u_2)\bigr) \sim\cases{ u_2^{-2/3},
&\quad$d=1$,
\vspace*{2pt}\cr
2u_2^{-1/2} \log^{1/2}(1/u_2),
&\quad$d = 2$,
\vspace*{2pt}\cr
u_2^{-1/2}, &\quad$d \ge3$,} %
\]
where we have written $a(n)$ for $a_n$ for readability.
Thus for (B2) we need $a(1/h_d(u_2)) \ll1/Nu_1 h_d(u_2)$, which means
$N \ll g_d(u_2)/u_1$.

The next order of business is to compute $\nu_d$.
Stochastic calculus (or calculations with infinitesimal generators)
tells us that
\[
v(x) = E_x \exp\biggl(-\int_0^{T_0}
Y_s \,ds \biggr) %
\]
is the unique function on $[0,\infty)$ with values in $[0,1]$,
$v(0)=1$ and
\[
v'' - x v = 0\qquad\mbox{in $d=1$}, \qquad
\beta_d x v'' - x v = 0 \qquad\mbox{in $d
\ge2$}. %
\]
In $d=1$ all solutions have the form:
\[
v(x) = \alpha Ai(x) + \beta Bi(x), %
\]
where $Ai$ and $Bi$ are Airy functions
\begin{eqnarray*}
Ai(x) & =& \frac{1}{\pi} \int_0^\infty\cos
\biggl( \frac{t^3}{3} + x t \biggr) \,dt,
\\
Bi(x) & =& \frac{1}{\pi} \int_0^\infty\exp
\biggl(- \frac{t^3}{3} + x t \biggr) + \sin\biggl( \frac{t^3}{3} + x t
\biggr) \,dt.
\end{eqnarray*}
Since $Bi$ is unbounded and $Ai$ is decreasing on $[0,\infty)$, we
take $\beta= 0$ and set $\alpha= 3^{2/3}\Gamma(2/3)$ to satisfy the
boundary condition, $v(0)=1$. Letting $\varepsilon\to0$ we conclude that
%
\begin{equation}
\gamma_1 = - \alpha Ai'(0)= 3^{1/3}
\Gamma(2/3)/\Gamma(1/3). \label{gam1}
\end{equation}
In $d \ge2$, $v(x) = \exp(-\beta_d^{-1/2} x)$, and we have
%
\begin{equation}
\gamma_d = \beta_d^{-1/2}. \label{gamd}
\end{equation}

\section{Proof, part II: Missing details for \texorpdfstring{$\lambda=1$}{lambda=1}}\label{sec4}

In the previous section we have calculated the probability $\nu
_d^\varepsilon$ that a type-1 family reaches size $\varepsilon/h_d(u_2)$
and then gives rise to a type 2. To let $\varepsilon\to0$ and prove
Theorem~\ref{th1} we need to consider the possibility of a mutation to
type 2 in a family that (i) never reaches size $n \varepsilon$, or (ii)
will reach $n \varepsilon$ but has not yet. To have a convenient name we
will call these small families. Families of the first kind arise at
rate $Nu_1(1-1/n\varepsilon)$ and families of the second kind arise at
rate $Nu_1/n\varepsilon$. We will now calculate the expected rate at
which type 2's are born from these small families. In the proof of
Theorem~\ref{th1}, we will let $\varepsilon\to0$ slowly as $n\to
\infty$ so we can and will assume $n\varepsilon\to\infty$.

Consider the voter model $\xi^0_t$ starting from a single 1 at the
origin at time 0. Let $V_k$ be the total time spent at level $k$, that is,
$|\{ t\dvtx  |\xi^0_t|=k\}|$ and let $N_k$ be the total number of returns
to level $k$ before leaving the interval $(0,n\varepsilon)$. Recalling
our assumption in (\ref{bdryA}), we let $q(k)$ the rate jumps occur at
level~$k$.

Let $S_k$ be the embedded discrete time chain, which is a simple random
walk, and let $T^+_k = \min\{ n \ge1\dvtx  S_n = k \}$.
%
\begin{eqnarray}\label{f1mh}
 E_1 \biggl( \int_0^{T_0} \bigl|
\xi^0_s\bigr| \,ds \bigg| T_0 < T_{n\varepsilon}
\biggr) &=& E_1 \Biggl(\sum_{k=1}^{n\varepsilon}
k V_k \bigg| T_0 < T_{n\varepsilon} \Biggr)\nonumber
\\
& =& E_1 \Biggl(\sum_{k=1}^{n\varepsilon}
\frac{k N_k}{q(k)} \bigg| T_0 < T_{n\varepsilon} \Biggr)
\\
& =&  \sum_{k=1}^{n\varepsilon}\frac{\overline P_1(T_k < \infty)}{\overline P_k(T_k^+>T_0)}
\frac{k}{q(k)},\nonumber
\end{eqnarray}
where the bar indicates conditioning on $T_0 < T_{n\varepsilon}$. A
similar argument shows that
%
\begin{equation}
E_1 \biggl( \int_0^{T_{n\varepsilon}} \bigl|
\xi^0_s\bigr| \,ds \bigg| T_{n\varepsilon} < T_0
\biggr) = \sum_{k=1}^{n\varepsilon}\frac{1}{\widehat P_k(T_k^+>T_{n\varepsilon
})}
\frac{k}{q(k)}, \label{f2mh}
\end{equation}
where the hat indicates conditioning on $T_{n\varepsilon}<T_0$.

The three conditional probabilities we need can be computed using facts
about simple random walk that follow from the fact that
it is a martingale.
%
\begin{equation}
\overline P_1(T_k<\infty) = \frac{P_1(T_k<\infty)P_k(T_0<T_{n\varepsilon
})}{P_1(T_0<T_{n \varepsilon})} =
\frac{(1/k)(1-k/n\varepsilon)}{(1-1/n\varepsilon)}. \label{cpr1}
\end{equation}
For the next two we note that the first step has to be in the correct
direction for these events to happen.
%
\begin{eqnarray}
\overline P_k\bigl(T_k^+>T_0\bigr) & =&
\frac{(1/2)(1/k)}{(1-k/n\varepsilon)}, \label{cpr2}
\\
\widehat P_k\bigl(T_k^+>T_{n\varepsilon}\bigr) & =&
\frac{(1/2)(1/({n\varepsilon-
k}))}{(k/n\varepsilon)}. \label{cpr3}
\end{eqnarray}
Thus the expected total man-hours $\int_0^{T_0} |\xi^0_s| \,ds$ for
a family that will die out before reaching size $n\varepsilon$ is
%
\begin{equation}
\label{PDieOut} \sim\frac{2}{(1-1/n\varepsilon)} \sum_{k=1}^{n\varepsilon}(1-k/n
\varepsilon)^2 \frac{k}{q(k)},
\end{equation}
and in families that have yet to reach size $n\varepsilon$,
%
\begin{equation}
\label{PNotYet} \frac{2}{n\varepsilon} \sum_{k=1}^{n\varepsilon}(n
\varepsilon-k) \frac
{k^2}{q(k)}.
\end{equation}

The next result shows that the contribution of small families are
indeed negligible. Note that in all three cases the order of magnitude
of the contributions from small families is the same as the overall
rate, but contains a constant that $\to0$ as $\varepsilon\to0$.

\begin{lemma} The expected total man-hours in small families is
\[
\le\cases{\displaystyle Nu_1u_2^{1/3} \cdot
\frac{\varepsilon^2}{4}, &\quad$d=1$,
\vspace*{2pt}\cr
\displaystyle Nu_1
u_2^{1/2} \log^{1/2} (1/u_2) \cdot
\frac
{7\varepsilon}{24\beta_2}, &\quad$d=2$,
\vspace*{2pt}\cr
\displaystyle Nu_1u_2^{1/2}
\cdot\frac{\varepsilon}{2d\beta_d}, &\quad$d\ge3$.} %
\]
\end{lemma}

\begin{pf}
In one dimension, $q(k) = 2$. The sum in (\ref{PDieOut}) is dominated by
\[
\int_0^{n\varepsilon}(1-x/n\varepsilon)^2x \,dx =
\frac{1}{(n\varepsilon)^2}\int_0^{n\varepsilon} y^2 (n
\varepsilon-y) \,dy = \frac{(n\varepsilon)^2}{12}. %
\]
Thus, families of the first kind produce type 2's at rate $\le Nu_1u_2
(n\varepsilon)^2/12$. The expression in (\ref{PNotYet}) is dominated by
\[
\frac{2}{n\varepsilon} \int_0^{n\varepsilon}(n\varepsilon-
x)x^2 \,dx = \frac{(n\varepsilon)^3}{6}. %
\]
Thus, families of the second kind produce type 2's at rate $\le Nu_1u_2
(n\varepsilon)^2/6$. Adding the last two conclusions gives the result for $d=1$.

In $d\ge3$, (\ref{bdryA}) implies $q(k) = 2d\beta_d k$, so (\ref
{PDieOut}) becomes
\[
\frac{1}{d \beta_d} \sum_{k=1}^{n\varepsilon}(1 -k/n
\varepsilon)^2. %
\]
The sum is bounded above by the integral
\[
\int_0^{n\varepsilon}(1-x/n\varepsilon)^2 \,dx =
\frac{n\varepsilon}{3}, %
\]
so with our choice of $n = u_2^{-1/2}$, families of the first kind
produce type 2's at rate bounded above by $Nu_1u_2^{1/2}\varepsilon
/(3d\beta_d)$.
Setting $q(k) = 2d\beta_d k$, (\ref{PNotYet}) becomes
\[
\frac{1}{d\beta_d n \varepsilon} \sum_{k = 1}^{n\varepsilon}(n
\varepsilon- k) k. %
\]
The sum is bounded above by the integral
\[
\int_0^{n\varepsilon}(n\varepsilon- x)x \,dx =
\frac{(n\varepsilon)^3}{6}. %
\]
Thus, families of the second kind produce type 2's at rate $\le
Nu_1u_2^{1/2}\varepsilon/(6d\beta_d)$. Adding the last two conclusions
gives the result for $d\ge3$.

In $d=2$, (\ref{bdryA}) implies $q(k) = 4 \beta_2 k/\log k$, so (\ref
{PDieOut}) becomes
\[
\frac{1}{2\beta_2}\sum_{k=1}^{n\varepsilon}(1-k/n
\varepsilon)^2\log k. %
\]
Each term in the sum is bounded above by $\log(n\varepsilon)$, so the
sum is less than $n\varepsilon\log n\varepsilon$. Since $n = u_2^{-1/2}\log
^{-1/2}(1/u_2)$, families of the first kind produce type 2's at rate
bounded above by
\begin{eqnarray*}
Nu_1u_2 \cdot\frac{1}{2\beta_2} n\varepsilon\log(n
\varepsilon) & =& Nu_1u_2 \cdot\frac{1}{2\beta_2} \varepsilon
u_2^{-1/2}\log^{-1/2}(1/u_2) \cdot
\frac{1}{2} \log(1/u_2)
\\
& =& \frac{\varepsilon}{4\beta_2} N u_1 u_2^{1/2}
\log^{1/2}(1/u_2).
\end{eqnarray*}
Taking $q(k) = 4 \beta_2 k/\log k$, (\ref{PNotYet}) becomes
\[
\frac{1}{2\beta_2n\varepsilon} \Biggl(\sum_{k=1}^{n\varepsilon}(n
\varepsilon- k)k \log k \Biggr). %
\]
The sum is bounded above by
\[
\int_0^{n\varepsilon} (n\varepsilon-x) x\log(n\varepsilon) \,dx
\le\frac{(n\varepsilon)^3 }{6}\log( n\varepsilon). %
\]
Thus families of the second kind produce type 2's at rate bounded above by
\begin{eqnarray*}
\frac{Nu_1u_2}{n\varepsilon} \cdot\frac{1}{2\beta_2n\varepsilon} \cdot\frac
{(n\varepsilon)^3 }{6}\log( n
\varepsilon) & =& \frac{1}{12\beta_2} Nu_1u_2 \cdot n\varepsilon
\log(n\varepsilon)
\\
& =& \frac{\varepsilon}{24\beta_2} N u_1 u_2^{1/2}
\log^{1/2}(1/u_2).
\end{eqnarray*}
Adding the last two conclusions gives the result for $d=2$ and
completes the proof.
\end{pf}

\section{Proof, part III: Almost neutral mutations}\label{sec5}

In the biased voter model, whose law we denote by $P^\lambda$, jumps
occur at rate $1+\lambda$ times the size of the boundary.
To compensate for this we need to run the unbiased ($\lambda=1$) voter
at rate $(1 +\lambda)/2$. If we do this,
call the resulting law $\widetilde P^0$, and let $\omega_T$ is a realization
of $\xi^0_t$ run up to time $T$ then the Radon--Nikodym derivative
\[
\frac{dP^\lambda}{d\widetilde P^0}(\omega_T) = \biggl( \frac{ 2\lambda
}{\lambda+ 1 }
\biggr)^{n_+} \biggl( \frac{ 2 }{\lambda+ 1 } \biggr)^{n_-},
\]
where $n_+$ and $n_-$ are the number of up jumps in $\omega_t$ when
$0\le t \le T$.

If $\max_{t \le T} |\xi^0_t| =O(K)$ then the difference $0 \le
n_+-n_- = O(K)$.
Since under $\widetilde P_0$, $|\xi^0_t|$ is a time change of simple
random walk, we see that
the total number of jumps $n_+ + n_- = O(K^2)$.
Taking $K=1/h_d(u_2)$ and assuming $|\lambda- 1| \ll h_d(u_2)$, when
$u_2$ is small the Radon--Nikodym derivative is
\begin{eqnarray*}
& =& \biggl( 1 + \frac{ \lambda- 1 }{\lambda+ 1 } \biggr)^{n_+} \biggl(
1 -
\frac{ \lambda- 1 }{\lambda+ 1 } \biggr)^{n_-}
\\
& =& \biggl( 1 + \frac{ \lambda- 1 }{\lambda+ 1 } \biggr)^{n_+-n_-}
\biggl( 1 -
\frac{ (\lambda- 1)^2 }{(\lambda+ 1)^2 } \biggr)^{n_+ +
n_-} \approx1.
\end{eqnarray*}

The last result implies that (\ref{procconv}) extends to almost
neutral mutations, and the computations in
Section~\ref{sec2} are valid. To extend the part of the proof in Section~\ref{sec3}, we
need to check that (\ref{cpr1})--(\ref{cpr3})
are true asymptotically for almost neutral mutations. To do this we
recall that if $a < x < b$
%
\begin{equation}
P^\lambda_x (T_b < T_a ) =
\frac{\theta^x - \theta^a}{\theta^b -
\theta^a} \qquad\mbox{where }\theta= 1/\lambda. \label{pqhit}
\end{equation}
When $0\le a < x \le b = O(1/h_d(u_2))$ and $|\lambda- 1| \ll
h_d(u_2)$ we have
\[
P^\lambda_x (T_b < T_a ) \approx
\frac{x-a}{b-a}. %
\]
To show that the sums come out the same we need the following uniform
version which follows from (\ref{pqhit}).
If $|\lambda- 1| h_d(u_2) \to0$ then for any $C$ fixed
\[
\sup_{0 \le-a, b \le C/h_d(u_2)} \biggl\llvert\frac{P^\lambda_0(T_b <
T_a)}{-a/(b-a)} - 1 \biggr\rrvert
\to0. %
\]



%

\printaddresses

\end{document}